\documentclass[11pt]{article}
\usepackage{amsmath,amssymb,amsfonts}
\usepackage[pdftitle={Title of Article}, verbose,
            colorlinks=true,
            naturalnames=true,
            linkcolor=blue,
            citecolor=red,
            pdfstartview=XYZ]{hyperref}
\newcommand{\Hom}{\mathrm{Hom}}

\newcommand{\Ext}{\mathrm{Ext}}

\newcommand{\Spec}{\mathrm{Spec}}

\newcommand{\Supp}{\mathrm{Supp}}

\newcommand{\depth}{\mathrm{depth}}
\newcommand{\width}{\mbox{width}\,}
\renewcommand{\dim}{\mathrm{dim}}

\newcommand{\pd}{\mathrm{pd}}
\newcommand{\id}{\mathrm{id}}

\newcommand{\Gid}{\mathrm{Gid}}

\newcommand{\E}{\mathrm{E}}

\newcommand{\T}{\mathrm}

\newcommand{\lo}{\longrightarrow}

\newtheorem{thm}{Theorem}[section]
\newtheorem{cor}[thm]{Corollary}
\newtheorem{lem}[thm]{Lemma}
\newtheorem{prop}[thm]{Proposition}
\newtheorem{defn}[thm]{Definition}

\numberwithin{equation}{section}
\begin{document} 
\title{Gorenstein injective dimension, Bass formula and
Gorenstein rings \footnotetext{The first author's research was partially 
supported by IPM grant No. 82130025.} \footnotetext {The second 
author's research was supported 
by IPM grant No. 82130212.} \footnotetext{{\it E-mail addresses:} 
\href{mailto:lkhatami@ictp.trieste.it}{lkhatami@ictp.trieste.it} and
\href{mailto:yassemi@ipm.ir}{yassemi@ipm.ir}}
} 
\author{Leila Khatami $^{\it a}$
and Siamak Yassemi 
$^{\it b,c}$
\\
{\small\it $^a$ The Abdus Salam ICTP, Strada Costiera 11, 34100
Trieste, Italy}\\
{\small\it $^b$ School of Mathematics, IPM, P.O. Box 19395-5746, Tehran, 
Iran}\\ 
{\small\it $^c$ Department of Mathematics, University of Tehran, P.O. Box 
13145-448, Tehran, Iran}
} \maketitle 
%\bigskip
%\vspace{2cm} 
\begin{abstract} \noindent Let $(R,\frak m, k)$ be a
noetherian local ring. It is well-known that $R$ is regular if and only if
the injective dimension of $k$ is finite. In this paper it is shown that
$R$ is Gorenstein if and only if the Gorenstein injective dimension of $k$
is finite. On the other hand a generalized version of the so-called Bass
formula is proved for finitely generated modules of finite Gorenstein
injective dimension. It also improves the results by Enochs and Jenda
\cite{EJbass} and Christensen \cite{ch}. \\
\\
{\it MSC:} 13D05; 13H10 
\end{abstract}
%\footnotetext{{\it Key words and phrases}:
%Gorenstein ring, Gorenstein injective dimension.}}
\section{Introduction} 
\vspace{.2in} 
The classical homological dimensions are no debt the central homological
notions in commutative algebra. The notion of Gorenstein injective
dimension of a module has been defined by E. E. Enochs and O. M. G. Jenda
\cite{EJ2} in mid nineties. It is a refinement of the classical injective
dimension and shares some of its nice properties. One can also consider
the Gorenstein injective dimension as the dual notion to the Gorenstein
dimension introduced by M. Auslander \cite{A}. In this note we try to
generalize some of the classical results on injective dimensions to
Gorenstein injective dimension.\\
Recall that the following statement is well-known.
\\
\\ 
\noindent{\bf Theorem} {\it If $(R,\frak m, k)$ is a commutative local
noetherian ring then $R$ is regular if and only if the injective dimension
of $k$ is finite.}
\bigskip 
\\ 
Enochs and Jenda poved that over a Gorenstein local
ring, the Gorenstein injective dimension of every module is finite
\cite{EJgor}. Using the so-called Foxby duality, it has been proved that
if the local ring $R$ admits a dualizing complex (i.e. it is a homomorphic
image of a Gorenstein local ring), then $R$ is Gorenstein if and only if
its residue field has finite Gorenstein injective dimension (cf. \cite{ch}
and \cite{chf}). In section 2, we prove the same statement over an
arbitrary noetherian local ring (Theorem \ref{t:gid}).\\ 
\\ 
The main theorem of section 3, generalizes the so-called Bass formula. 
Recall that\\
\\
\noindent{\bf Theorem} {\it If  $(R, \frak{m}, k)$ be a
commutative noetherian local ring and $M$ is a finitely generated
$R$-module of finite injective dimension, then $$\id _R
M=\depth_R M.$$} 

In (\cite{ch}; 6.2.15) Christensen has proved the
same formula for finitely generated modules of finite Gorenstein
injective dimension over a Cohen-Macaulay local ring which admits
a dualizing module. More recentely, the result has been proved over 
local rings which admit a dualizing module (cf. \cite{chf}). 
\\
\noindent Theorem \ref{t:gbass} gives another generalization of the Bass
formula. Namely, \\
\\
\noindent{\bf Theorem} {\it Let $S$ be a commutative 
noetherian ring. If $M$ is a finitely generated
$S$-module of finite Gorenstein injective dimension, then} $$\Gid_S M =
\sup \{ \, \depth \, S_{\frak p} \, | \, \frak{p} \in \Supp(M) \, \}.$$

\noindent Cosequently, the equation $\Gid_RM =\depth \, R$ holds, for a
finitely generated module $M$ of finite Gorenstein injective dimension 
over an almost Cohen-Macaulay local ring $R$ (\ref{cor:cmd}).\\
\\
\noindent{\bf {Convention.}} Throughout this note, the rings are assumed
to be commutative and notherian. Furthermore, $(R,\frak m, k)$, always
denotes a commutative noetherian local ring with the maximal ideal $\frak 
m$ and the residue field $k$.

%\section{Basic Definitions}
%\vspace{.2 in}
\section{Characterization of Gorenstein local rings.}
\vspace{.2in} In this section we prove that finiteness of the
Gorenstein injective dimension characterizes the Gorenstein local
rings. First recall basic definitions and facts. For details and
proofs cf. \cite{EJr} or  \cite{ch}.
\begin{defn}
An exact complex of injective $R$-modules,
$$I=\ldots \to I_2\lo I_1\lo I_0\lo I_{-1}\lo
I_{-2}\lo\ldots$$ is a {\it complete injective resolution} if and
only if the complex $\T{Hom}_R(J,I)$ is exact for every injective
$R$-module $J$. A module $M$ is said to be {\rm Gorenstein
injective} if and only if it is the 0-th kernel of a complete
injective resolution.
\end{defn}
It is clear that every injective module is Gorenstein injective.
Then one can construct a Gorenstein injective resolution for any
module.

\begin{defn}
Let $M$ be an $R$-module. A {\rm Gorenstein injective resolution}
of $M$ is an exact sequence $$0 \to M \to G_0 \to G_1 \to
\cdots$$ such that $G_i$ is Gorenstein injective for all $i \geq
0$. We say that $M$ has {\rm Gorenstein injective dimension} less
than or equal to $n$, $\Gid_R M \leq n$, if $M$ has a Gorenstein
injective resolution $$0 \to M \to G_0 \to G_1 \to \cdots \to G_n
\to 0.$$
\end{defn}
It is well-known that one always has $$\Gid_R M \leq \id_R M.$$
The equality holds if $\id_RM< \infty$. The Gorenstein injective
dimension can be computed using $\Ext$ functors.
\begin{thm}
Let $M$ be an $R$-module of finite Gorenstein injective
dimension. Then
$$\Gid_RM = \sup \{i\, | \,  \Ext_R^i(J,M) \neq 0 \,
\, \mathrm{for \, \, an}  \, \, R\mathrm{-module}\, \, J \, \,
\mathrm{with}\, \, \id_RJ<\infty \}.$$
\end{thm}
We start with proving some preliminary lemmas. \\ Recall that over a 
local noetherian ring, $(R, \frak m,k)$, injective dimension of a finite 
$R$-module $M$ is the supremum of integers $i$ such that $\Ext_R^i(k,M)$ 
is non-zero (cf. \cite{BH}). The following lemma shows that the residue
field can be 
replaced with its injective envelope when the module has finite injective 
dimension.
\begin{lem}\label{lemid} Let $\phi : (R,\frak m,k)\rightarrow (S, \frak
n, l)$ be a local ring homomorphism of noetherian local rings and let $M$
be a finitely generated $S$-module with finite injective dimension over
$R$. Then $$\id_R M=\sup\{i\, | \,\Ext^i_R(\E(k),M)\neq 0\}.$$ 
\end{lem}
\noindent{\it Proof.} Let $\id_R M=t$. The exact sequence $$0\to
k\to \E(k)\to C\to 0$$ induces the long exact sequence
$$\cdots
\to\Ext^t_R(\E(k),M)\to\Ext^t_R(k,M)\to\Ext^{t+1}_R(C,M)\to
\cdots.$$ Since $\Ext^{t+1}_R(C,M)=0$ and $\Ext^t_R(k,M)\neq 0$
(cf. \cite{AF}; 5.5) one has $$\Ext^t_R(\E(k),M)\neq 0$$ and this
proves the assertion.\hfill$\square$

\bigskip
\noindent Note that in \ref{lemid} the finitely generated
condition for $M$ is necessary. \\
\\ 
\noindent{\bf Example.} Let $\phi$ be the identity
homomorphism over $R$. If $M=\E(R/\frak p)$, then $\id_R M=0$. Let
$\psi:\E(k)\to\E(R/\frak p)$ be an $R$ homomorphism. If $x\in\E(k)$ then
module $Rx$ is of finite length. Thus $R\psi(x)$ is a submodule of
$\E(R/\frak p)$ which has finite length . Since $\E(R/\frak p)$ is an
essential extension of $R/\frak p$ non of its non-trivial submodules have
finite length. Therefore $\psi(x)=0$. That is $\Hom_R(\E(k),\E(R/\frak
p))=0$.

\begin{lem}\label{lem}
Let $(R,\frak m,k)$ be a noetherian local ring and let $M$ be an
$R$-module. Then for any $\frak p\neq\frak m$ and any $i\ge 0$,
$\Ext^i_R(\E (R/ {\frak p}),M)=0$ if one of the following
conditions hold.
\begin{itemize}
\item[(a)] $M$ has finite length.
\item[(b)] $R$ is complete and $M$ is finitely
generated.
\end{itemize}
\end{lem}
\noindent{\it Proof.} (a) Since $\frak p\neq\frak m$ there exists
$x\in R-\frak p$ such that $xM=0$. But multiplication by $x$ is an
automorphism on $\E(R/\frak p)$. So the assertion holds.\\
(b) See the proof of (\cite{FFGR}; 2.2). \hfill$\square$\\
\\
The following corollary, which genralizes \ref{lemid}, shows that the
Gorenstein injective dimension of a finith-lengh module, if it is finite,
can be computed in terms of vanishing of the $\Ext_R^i(\E(k),_)$
functors.
\\
\begin{cor}\label{c:ee}
Let $(R,\frak m,k)$ be a noetherian local ring and let $M$ be an
$R$-module of finite length which has finite Gorenstein injective
dimension. Then
$$\Gid_R M=\sup\{i\, | \, \Ext^i_R(\E(k),M)\neq 0\}.$$
\end{cor}
\noindent{\it Proof.} It is proved (\cite{H}; 2.29) that for any
$R$-module $M$ with finite Gorenstein injective dimension we have
$$\Gid_R M=\sup\{i\, | \, \exists \frak p\in\Spec R:
\Ext^i_R(\E(R/\frak p),M)\neq 0\}. $$ Now the assertion follows
from \ref{lem}.\hfill$\square$
\\
\\
Now we are ready to prove the main theorem of this section. It gives a
characterization of Gorenstein local rings in terms of the finiteness of
Gorenstein dimension of modules. In \cite{chf}, using hyperhomological
techniques, a similar result has been proved, provided that the ring
admits
a dualizing complex (equivalently, is a homomorphic image of a Gorenstein
local ring). Theorem \ref{t:gid} generalizes that theorem with a rather
simpler proof.
\begin{thm}\label{t:gid}
Let $(R,\frak m,k)$ be a noetherian local ring. The following are
equivalent.
\begin{itemize}
\item[(i)] $R$ is Gorenstein.
\item[(ii)] $\Gid_R k$ is finite.
\item[(iii)] $\Gid_R M$ is finite for any finitely
generated $R$--module $M$.
\item[(iv)] $\Gid_R M$ is finite for any $R$-module
$M$.
\end{itemize}
\end{thm}
\noindent{\it Proof.} (iv)$\Rightarrow$(iii) and
(iii)$\Rightarrow$(ii) are obvious. For (ii)$\Rightarrow$(i), let
$i\ge 0$. One has
\[ \begin{array}{rl} \Ext^i_R(\E(k), k)
& \, =\Ext^i_R(\E(k),\Hom_R(k,\E(k)))\\
& \, \cong\Ext^i_R(k,\Hom_R(\E(k),\E(k)))\\
& \, \cong\Ext^i_R(k, \hat{R})\\
& \, \cong\Ext^i_{\hat{R}}(k,\hat{R})
\end{array} \]
where $\hat{R}$ is the completion of $R$ in $\frak m$-adic
topology. Now using \ref{c:ee}, we have
$$\id_{\hat{R}}
\hat{R}=\Gid_R k < \infty.$$ Thus $\hat{R}$ and hence $R$ are
Gorenstein.\\ (i)$\Rightarrow$(iv). cf. \cite{EJgor} or \cite{ch}.
\hfill$\square$
\\
\\
Let $R$ be a local ring. It is well-known that $R$ is {\it regular} 
(respectively, {\it Gorenstein}, {\it Cohen-Macaulay}) if and only if 
there exists a {\it simple} (respectively, {\it cyclic}, {\it finitely 
generated}) $R$--module with finite injective dimension. (For the first
and last statement cf. \cite{BH} and for the second one cf. \cite{PS}.)
\\
Theorem \ref{t:gid} shows that if there exists a simple $R$-module of 
finite
Gorenstein injective dimension then $R$ is Gorenstein. Now it is
natural to ask
\\
\\
\noindent{\bf Question.} What one can say about a ring $R$ which
admits a cyclic (respectively, finitely generated) module of finite
Gorenstein injective dimension?
\section{Bass formula.}
Recall that if a finitely generated module over a noetherian
local ring has finite injective dimension then its injective
dimension is equal to the depth of the base ring. This is known
as Bass formula. In (\cite{ch}; 6.2.5) Christensen has proved
that over a Cohen-Macaulay local ring with a dualizing module, one
can replace injective dimension with Gorenstein injective
dimension. In this section we try to generalize this result. The
main result of this section is the following theorem.

\begin{thm}\label{t:gbass}
Let $S$ be a commutative noetherian ring. If $M$ is a finitely 
generated
$S$-module of finite Gorenstein injective dimension, then $$\Gid_S M =
\sup \{ \, \depth \, S_{\frak p} \, | \, \frak{p} \in \Supp(M) \, \}.$$
\end{thm}

\noindent{\it Proof.} Let $\Gid_SM=0$ and suppose that $\depth \, S_{\frak 
p} >0$ for some $\frak {p} \in \Supp(M)$. 
Then there exists an 
$S_{\frak p}$-regular
element $x \in \frak{p}S_{\frak p}$. The exact sequence $$0 \to S_{\frak 
p} 
\stackrel{.x}{\to} 
S_{\frak p} \to S_{\frak p}/xS_{\frak p} \to 0$$ induces the following 
exact
sequence.$$ M_{\frak p} \stackrel{.x}{\to} M_{\frak p} \to 
\Ext_{S_{\frak p}}^1(S_{\frak p}/xS_{\frak p} , M_{\frak p}) \to 
0$$ 
Since $M$ is finitely generated, using Nakayama's
lemma, multiplication by $x$ is not surjective over $M_{\frak p}$ and so
$\Ext_{S_{\frak p}}^1(S_{\frak p}/xS_{\frak p} , M_{\frak p}) \neq 0$. \\
On the other hand, $M$ is a Gorenstein injective $R$-module and then an 
exact sequence $0 \to N \to I \to M \to 0$, with $I$ an
injective $S$-module, exists. Hence $\Ext_{S_{\frak p}}^1(S_{\frak 
p}/xS_{\frak p}, M_{\frak p}) \cong 
\Ext_{S_{\frak p}}^2(S_{\frak p}/xS_{\frak p}, N_{\frak p}) = 0$, which is 
a contradiction. Then $\depth\,
S_{\frak p}$ should be zero and the claim is proved in this case.

Now let $\Gid_SM=n>0$. We prove that $M$ can be written as a homomorphic
image of an $S$-module with injective dimension $n$. Furthermore, we will
show that localizations of this module to every prime ideal outside 
$\Supp(M)$ is injective.
\\
By (\cite{H}, 2.45), there exists a 
Gorenstein injective $S$-module $G$ and an $R$-module $C$ with $\id_S C
=\Gid_S C = n-1$, such that the following sequence is exact.$$ 0 
\rightarrow M \rightarrow G \rightarrow C
\rightarrow 0$$
Since $G$ is a Gorenstein injective $S$-module,
there exists an injective $S$-module $E$ and an exact sequence
$$0 \rightarrow K \rightarrow E \rightarrow G
\rightarrow 0 $$ with $K$ a Gorenstein injective $S$-module, too.
So the isomorphisms $C \cong G/M$ and $G
\cong E/K$ of $S$-modules hold and then there exists a submodule $L$ of 
$E$ such that $K
\subseteq L$ and $M \cong L/K$ and therefore $C \cong E/L$.

Considering the following exact sequence, we get $\id_S L 
\leq \id_S C+1=n$.
$$0 \rightarrow L \rightarrow E \rightarrow C
\rightarrow 0 $$

On the other hand, since $\Gid_S M =n$, there exists an injective
$S$-module $J$ with $\T{Ext}_S^n(J,M)\neq 0$.
The exact sequence $$0 \rightarrow K \rightarrow L \rightarrow M
\rightarrow 0$$ induces the following exact sequence.
$$\T{Ext}_S^n(J,L) \rightarrow \T{Ext}_S^n(J,M)
\rightarrow \T{Ext}_S^{n+1}(J,K)=0$$ Therefore, we have 
$\T{Ext}_S^n(J,L)\neq
0$ and then $\id_S L \geq n$. Hence $\id_S L=n=\Gid_SM$.
By Chouinard's equality (\cite{C}), we have $$\id_SL\, =\,\sup\,  \{ 
\, \depth \, S_{\frak p}-\width_{S_{\frak
p}}L_{\frak p} \, | \, {\frak
p}\in \Supp(L)\, \}.$$
For every $\frak {p}\in \Supp(L)$, we have the following exact sequence. 
$$0\to K_{\frak p} \to 
L_{\frak p} \to M_{\frak p} \to 0$$ 
If $\frak {p} \in \Supp(M)$, then from the induced exact sequence 
$L_{\frak p}/{\frak p} L_{\frak p}\to M_{\frak p}/{\frak p} M_{\frak p}
\to 0$ and the fact that $M_{\frak p}/{\frak p}M_{\frak p} \neq 0$, 
we get $\width_{S_{\frak p}}L_{\frak p}=0$. \\
If ${\frak p} \not \in \Supp(M)$, then $L_{\frak p}\cong K_{\frak p}$. 
Therefore $L_{\frak p}$ is an injective $S_{\frak p}$-module and
then using Chouinard's equality, we have $$\depth \, S_{\frak 
p}-\width_{S_{\frak p}}L_{\frak p} \leq \id_{S_{\frak p}}L_{\frak p}=0.$$
Therefore
$$\begin{array}{ll}
\Gid_SM & = \id_SL \\
& =\sup\,  \{ \, \depth \, S_{\frak p}-\width_{S_{\frak p}}L_{\frak p} \, 
| \, {\frak p}\in \Supp(L)\, \} \\
& =\sup\,  \{ \, \depth \, S_{\frak p} \, | \, {\frak p}\in \Supp(M)\, \} 
\end{array}$$
This finishes the proof.
\hfill$\square$
\begin{cor}\label{cor:cmd}
Let $(R, \frak{m}, k)$ be an almost Cohen-Macaulay local ring (i.e. $\dim
\, R -\depth \, R \leq 1$) and let $M$ be a finite $R$-module. If $\Gid_RM
<\infty$ then $$\Gid_RM = \depth \, R.$$ 
\end{cor}
\noindent{\it Proof.} Use \ref{t:gbass} and the fact that over an almost 
Cohen-Macaulay ring, for every two prime ideals $ \frak{p} $ and 
$ \frak{q} $ with $ \frak{p} \in \frak{q} $, the inequality 
$\depth \, R_{\frak p} \leq \depth \, R_{\frak q}$ holds. \hfill$\square$
\\
\\
In \cite{SSS}, Salarian, Sather-Wagstaff and Yassemi prove that over local
rings, Gorenstein injective dimension "behaves well with respect to
killing a regular element". Namely, If $M$ is a finitely generate module
over a local ring $(R,\frak m, k)$ then $\Gid_R M < \infty$ implies
$\Gid_{R/xR} M/xM < \infty$, when $x\in \frak m$ is $R$- and $M$-regular.
The similar result was proved in \cite{chf} when $R$ is assumed to admit a
dualizing complex.
The next corollary of \ref{t:gbass} expresses how the values of the
Gorenstein injective
dimensions relate.
\begin{cor}
Let $(R,\frak m, k)$ be a noetherian local ring and $M$ a finitely
generated $R$-module. If $x\in \frak m$ is an $R$- and $M$-regular
element, then $$ \Gid_{R/xR} M/xM \leq \Gid_R M - 1.$$ Furthermore, the
equality holds when $R$ is almost Cohen-Macaulay and $\Gid_R M$ is finite. 
\end{cor}
\noindent{\it Proof.} If $\Gid_R M$ is not finite then the inequality is
clear. 
\noindent Now assume that $M$ has finite Gorenstein injective dimension,
then to
prove the
inequality, it is sufficient to use \ref{t:gbass} and the follwing facts.
\begin{itemize}
\item $\Supp(M/xM) = \{ \, {\frak p}/xR \, | \, {\frak p} \in
\Supp(M) \, \, \mathrm{and} \, \, x\in {\frak p} \, \}.$
\item If $x \in {\frak p}$ then $\depth \, (R/xR)_{{\frak p}/xR}
=
\depth \, R_{\frak p} -1.$
\end{itemize}  
The last part of the corollary is clear from \ref{cor:cmd}. 
\hfill$\square$ 
\\
\\
The following immediate corollary of \ref{t:gbass} shows that finite
Gorenstein injective dimension does not exceed after localization. 
\begin{cor} Let $S$ be a commutative noetherian ring and $M$ a finitely
generated $S$-module such that $\Gid_{S_{\frak p}} M_{\frak p} < \infty$
for all prime ideals $\frak p$. Then if $\frak p$ and $\frak q$ are two
prime ideals with $\frak p \subset \frak q$, we have $$\Gid_{S_{\frak p}}
M_{\frak p} \leq \Gid_{S_{\frak q}} M_{\frak q}.$$
\end{cor} 
\bigskip
\noindent{\bf Remark.} Let $(R,\frak m,k)$ be a noetherian local
ring and let $M$ and $N$ be finitely generated $R$-modules. If $\id_R
N<\infty$ then in \cite{Isc} Ischebeck showed that $$\depth \, R-\depth_R
M=\sup\{i\, | \,\Ext^i_R(M,N)\neq 0\}.$$ 
It is natural to ask whether one could replace $N$ by a finite module of
finite Gorenstein dimension. The answer is negative.
\\ 
\\
\noindent{\bf Example.} If $M$ has finite projective dimension then the
above equality holds for every finitely generated $R$-module $M$. Now let
$R$ be a Gorenstein local ring which is not regular. Let $k$ be the
residue field of $R$. One has $\pd_R k = \infty$ and then $ \sup\{i\, | \,
\Ext^i_R(k,k)\neq 0\}= \infty$. But it is clear (cf. \ref{t:gid}) that
$\Gid _R k < \infty$.
\\
The following
statement is a partial generalization of the Ischebeck's result
in another direction. It proves the equality for an $S$-module
$M$ (not necessarily finitely generated) with $\depth_R M=0$.
\vspace{.1in}
\begin{prop}\label{propis}
Let $\phi : (R,\frak m,k) \to (S, \frak n, l)$ be a local ring
homomorphism of noetherian local rings and let $M$ be an
$S$-module. Then for any finitely generated $S$-module $N$ with
finite injective dimension over $R$ the following equality holds.
$$\depth \, R - \depth_R M =\sup\{i\, | \,
\Ext_R^i(M,N)\neq 0\}$$ Provided that $\depth_R M=0$ or $M$ is
finitely generated $R$-module.
\end{prop}
\vspace{.1in} \noindent{\it Proof.} Set $\id_R N=t$. If $\depth_R
M=0$ then there exists a short exact sequence $$0\to k\to M\to
C\to 0$$ which induces a long exact sequence
$$\cdots\to\Ext^t_R(M,N)\to\Ext^t_R(k,N)\to\Ext^{t+1}_R(C,N)\to\cdots.$$
Since $\Ext^{t+1}_R(C,N)=0$ and $\Ext^t_R(k,N)\neq 0$, we have
 $\Ext^t_R(M,N)\neq 0$ and then $$ \sup\{i\, | \,
\Ext^i_R(M,N)\neq 0\}\ge t.$$ The inverse inequality holds
clearly. If $M$ is finitely generated we use induction on
$\depth_R M$ to prove the desired equality. If $\depth_RM >0$ then
there exist an $M$-regular element $x \in \frak m$. Using the long
exact sequence induced by the following exact sequence, the
equality can be proved from induction's hypothesis. $$0 \to M
\stackrel{.x}{\rightarrow} M  \to M/xM \to 0$$  \hfill$\square$

The following corollary of \ref{propis} is a generalization of
Bass formula. The statement has been appeared in \cite{TY}, too.
\begin{cor}
Let $\phi : (R,\frak m, k) \to (S, \frak n, l)$ be a local ring
homomorphism of noetherian local rings. For any finitely
generated $S$-module $N$ with finite injective dimension over
$R$, the following equality holds.$$\depth \, R=\id_R N$$
\end{cor}
\noindent{\it Proof.} In \ref{propis}, set $M=R/\frak{m}$ and use
(\cite{AF}; 5.5). \hfill$\square$
\bibliography{ref}
\nocite{BH,EJr,EJ1,EJ3,EJ0}
\bibliographystyle{plain}

%\parbox[t]{7cm}{Leila Khatami\\ The Abdus Salam ICTP \\
%Strada Costiera 11 \\ 34100 Trieste, Italy\\
%{\it E-mail:}
%\href{mailto:lkhatami@ictp.trieste.it}{lkhatami@ictp.trieste.it}
%}
%\hfill
%\parbox[t]{6.5cm}{Siamak Yassemi \\ Institute for Studies in Theoretical
%Physics and Mathematics (IPM) \\ P.O. Box 19395-5746, Tehran, Iran \\
%and \\
%Department of
%Mathematics \\ University of
%Tehran \\ P.O. Box 13145-448, Tehran, Iran \\
%{\it E-mail:} \href{mailto:yassemi@ipm.ir}{yassemi@ipm.ir}
%}
\end{document}